\documentclass[final,1p,times]{elsarticle}
\usepackage{bbm, amsmath,amssymb}
\usepackage{epstopdf}
\usepackage{graphicx}
\usepackage{psfrag}

\usepackage{mathrsfs}
\DeclareSymbolFontAlphabet{\mathrsfs}{rsfs}
\DeclareMathAlphabet{\mathcal}{OMS}{cmsy}{m}{n}

\newcommand{\be}{\begin{equation}}
\newcommand{\ee}{\end{equation}}

\journal{Journal of Computational Physics}

\begin{document}
\begin{frontmatter}

\title{Hyperboloidal layers for hyperbolic equations \\ on unbounded domains}
\author{An{\i}l Zengino\u{g}lu}
\address{Laboratoire Jacques-Louis Lions, Universit\'e Pierre et Marie Curie (Paris VI), Paris, France, and \\ Theoretical Astrophysics, California Institute of Technology, Pasadena, CA, USA}

\begin{abstract}
We show how to solve hyperbolic equations numerically on unbounded domains by compactification, thereby avoiding the introduction of an artificial outer boundary. The essential ingredient is a suitable transformation of the time coordinate in combination with spatial compactification. We construct a new layer method based on this idea, called the hyperboloidal layer. The method is demonstrated on numerical tests including the one dimensional Maxwell equations using finite differences and the three dimensional wave equation with and without nonlinear source terms using spectral techniques.
\end{abstract}

\begin{keyword}
Transparent (nonreflecting, absorbing) boundary conditions \sep perfectly matched layers \sep hyperboloidal layers \sep hyperboloidal compactification\sep wave equations \sep Maxwell equations. 
\end{keyword}

\end{frontmatter}

%%%%%%%%%%%%%%%%%%%%%%%%%%%%%%%%%%%

\section{Introduction}
\label{sec:intro}
Hyperbolic equations typically admit wavelike solutions that oscillate infinitely many times in an unbounded domain. Take a plane wave in one spatial dimension with frequency $\omega$ and wave number $k$,
\be\label{plane} u(x,t) = e^{2\pi i (k x - \omega t)}\,. \ee
Any mapping of such an oscillatory solution from an infinite domain to a finite domain results in infinitely many oscillations near the domain boundary, which can not be resolved numerically. We refer to this phenomenon as the compactification problem  \cite{GrOrs}. It is commonly stated that hyperbolic partial differential equations are not compatible with compactification, and therefore can not be solved on unbounded domains accurately. 

A suitable transformation of the \emph{time coordinate}, however, leads to a finite number of oscillations in an infinite spatial domain. Introduce 
\be \label{tau} \tau(x,t) = t - \frac{k}{\omega}\,\left(x + \frac{C}{1+x}\right)\,,\ee
where $C$ is a finite, positive constant. The plane wave \eqref{plane} becomes
\be\label{newplane} u(x,\tau) = e^{-2\pi i \left(k C/(1+x)+\omega \tau\right)}\,. \ee
This representation of the plane wave has only $k\, C$ cycles along a constant time hypersurface in the unbounded space $x\in[0,\infty)$, and is therefore compatible with compactification.

The simple idea just described has far reaching consequences. In numerical calculations of hyperbolic equations one typically truncates the unbounded solution domain by introducing an artificial outer boundary that is not part of the original problem. Boundary conditions, called transparent, absorbing, radiative, or nonreflecting, are constructed to simulate transparency of this artificial outer boundary. There has been significant developments in the treatment of artificial outer boundaries since the 70s, but there is no consensus on an optimal method \cite{Givoli, HagstromLau}. Especially the construction of boundary conditions for nonlinear problems is difficult \cite{Szeftel}. A successful technique for numerical calculations on unbounded domains resolves this problem for suitable hyperbolic equations and provides direct quantitative access to asymptotic properties of solutions.

Furthermore, the numerical construction of oscillatory solutions as \eqref{newplane} can be very efficient. Numerical accuracy requirements for hyperbolic equations are typically given in terms of numbers of grid points per wavelength. In the example presented above, the free parameter $C$ determines the number of cycles to be resolved, which can be chosen small. This suggests that high order numerical discretizations requiring a few points per wavelength can be very efficient in combination with time transformations of the type \eqref{tau}. 

The rest of the paper is devoted to the discussion of time transformation and compactification for hyperbolic equations. The theoretical part of the paper (sections \ref{sec:2} and \ref{sec:mult}) includes a detailed description of the method. We discuss the compactification problem (section \ref{sec:sp}) and its resolution (section \ref{sec:comp}) for the advection equation in one dimension. In section \ref{sec:wave} we discuss the wave equation with incoming and outgoing characteristics. We show that the method works also for systems of equations (section \ref{sec:sys}).  Hyperboloidal layers are introduced in section \ref{sec:layer} in analogy to absorbing layers. In multiple spatial dimensions, compactification is performed in the outgoing direction in combination with rescaling to take care of the asymptotic behavior (sections \ref{sec:resc} and \ref{sec:conf}). The layer strategy in multiple dimensions allows us to employ arbitrary coordinates in an inner domain, where sources or scatterers with irregular shapes may be present (section \ref{sec:mult_lay}). We finish the theoretical part discussing possible generalizations of the method to nonspherical coordinate systems (section \ref{sec:nonsph}). Section \ref{sec:num} includes numerical experiments in one and three spatial dimensions. In one dimension, we solve the Maxwell equations using finite difference methods (section \ref{sec:one}). A  stringent test of the method is the evolution of off-centered initial data for the wave equation in three spatial dimensions with and without nonlinear source terms (section \ref{sec:three}). We conclude with a discussion and an outlook in section \ref{sec:disc}.

\section{Compactification in one spatial dimension} \label{sec:2}
\subsection{Spatial compactification} \label{sec:sp}
\begin{figure}
\center
\includegraphics[height=0.17\textheight]{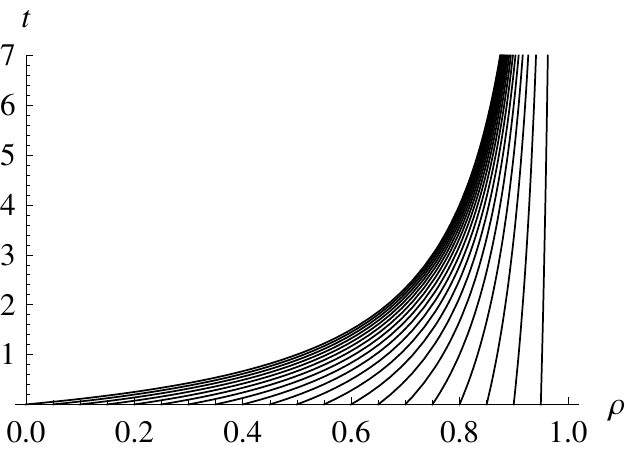}\hspace{1cm}
\caption{Characteristic diagram for the advection equation after spatial compactification. The characteristic speed approaches zero near spatial infinity at $\{\rho=1\}$ causing loss of numerical resolution: this is the compactification problem.\label{fig:1}}
\end{figure}

Consider the initial boundary value problem for the advection equation
\be\label{eq:advection}
\partial_t u + \partial_x u = 0, \qquad u(x,0) = u_0(x), \quad u(0,t) = b(t).  
\ee
The problem is posed on the unbounded domain $x\in[0,\infty)$. We transform the infinite domain in $x$ to a finite domain by introducing the compactifying coordinate $\rho$ via
\be\label{comp} \rho(x) = \frac{x}{1+x}, \qquad x(\rho) = \frac{\rho}{1-\rho}. \ee
%\be \label{eq:comp1} \rho: x \mapsto \rho(x) = \frac{-C+\sqrt{C^2+x^2}}{x} \quad \Rightarrow \quad x(\rho) = \frac{2\, C\, \rho}{1-\rho^2}. \ee
The advection equation becomes
\be \label{eq:adv} \partial_t u + (1-\rho)^2 \partial_\rho u = 0\,.\ee
%\be \label{eq:adv} \partial_t u = -\frac{(1-\rho^2)^2}{2(1+\rho^2)} \,\partial_\rho u, \quad \textrm{on} \quad \mathcal{D}= \{\rho\in[0,1),\ t\in[0,T)\}. \ee
The spatial domain is now given by $\rho\in[0,1]$. Characteristics of this equation are solutions to the ordinary differential equation
\[ \frac{d\rho(t)}{dt} = -(1-\rho(t))^2. \]
%\[ \frac{d\rho(t)}{dt} = -\frac{(1-\rho(t)^2)^2}{2(1+\rho(t)^2)}\,. \]
%For initial data $\rho(0)=\rho_0$ they read
%\[ \rho(t) =  \frac{\rho_0-t(1-\rho_0)}{1-t(1-\rho_0)}. \]
They are plotted in figure \ref{fig:1}. The compactification problem is clearly visible: the coordinate speed of characteristics approaches zero near a neighborhood of the point that corresponds to spatial infinity. The advection equation has a finite speed of propagation, therefore its characteristics can not reach infinity in finite time.

A concrete example illustrates the problem for oscillatory solutions. Set initial data $u_0(x) =\sin(2\pi x)$ and boundary data $b(t) = -\sin(2\pi t)$ in (\ref{eq:advection}). We obtain the solution
\be\label{eq:adv_sol} u(x,t) = \sin(2\pi(x-t)),\ee
which reads in the compactifying coordinate \eqref{comp}
\be\label{eq:adv_comp} u(\rho,t) = \sin \left(2\pi\left(\frac{\rho}{1-\rho} - t \right)\right). \ee
%\be\label{eq:adv_comp} u(\rho,t) = \sin\left(\frac{2\,\rho}{1-\rho^2} - t \right). \ee
\begin{figure}[ht]
\center
\includegraphics[height=0.17\textheight]{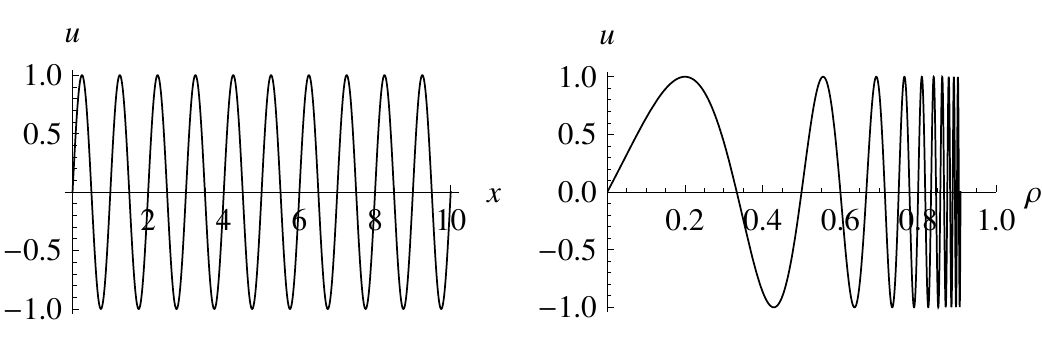}
\caption{The solution (\ref{eq:adv_sol}) at time $t=0$ is plotted on the left panel. The same solution in the compactifying coordinate given in (\ref{eq:adv_comp}) plotted on the right panel illustrates infinite blueshift in frequency.\label{fig:blue}}
\end{figure}
The solution is depicted in figure \ref{fig:blue} at $t=0$ on $x\in[0,10]$ in the original coordinate and on $\rho\in[0,10/11]$ in the compactifying coordinate.  The oscillations can not be resolved in the compactifying coordinate near infinity due to infinite blueshift in spatial frequency.

Mapping infinity to a finite coordinate seems to require infinite resolution. %A common way to deal with this problem is to add artificial viscosity to the equation so that spurious reflections near infinity are damped \cite{AppCol09}. 
However, a suitable time transformation discussed in the next section provides a clean solution to this problem.

\subsection{Hyperboloidal compactification}\label{sec:comp}
\begin{figure}[ht]
\center
\includegraphics[height=0.17\textheight]{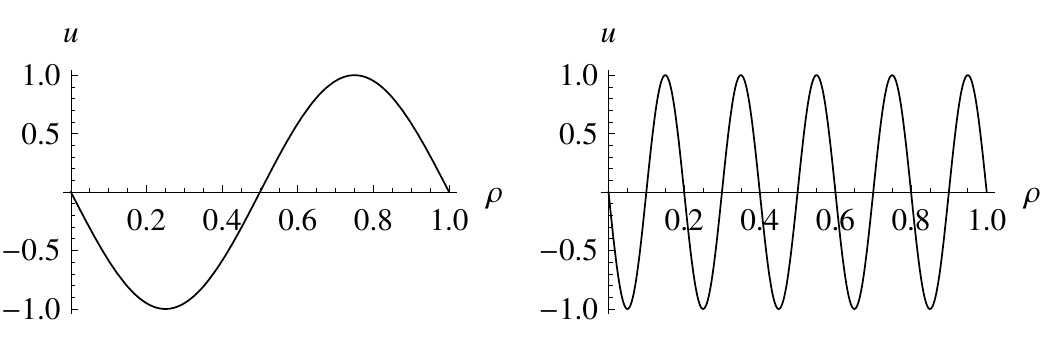}\hspace{1cm}
\caption{The solution \eqref{eq:hypsol} at time $\tau=1$ is plotted for two values of $C$. On the left panel we have $C=1$ and on the right panel $C=5$. The number of oscillations, and therefore the wavelength of the solution, can be influenced by the free parameter $C$.
\label{fig:red}}
\end{figure}

The idea is to transform the time coordinate as in \eqref{tau}. We introduce
\be\label{tauu} \tau = t - \left(x + \frac{C}{1+x}\right)\,,\ee
With the compactification \eqref{comp} we get the Jacobian
\[ \partial_\tau = \partial_t, \qquad \partial_x = (-1 + C \, \Omega^2) \partial_\tau + \Omega^2\,\partial_\rho\,,  \]
where we define $\Omega := 1-\rho$. The advection equation in the new coordinates $(\rho,\tau)$ reads
%\[ \partial_\tau u = - \sqrt{(1-\rho)^2+\rho^2} (\sqrt{(1-\rho)^2+\rho^2} - \rho) \partial_\rho u. \]
%\be\label{eq:hyp_adv} \partial_\tau u = - \frac{1}{2C} \,(1+\rho^2)\,\partial_\rho u.\ee
\[ \partial_\tau u + \frac{1}{C} \partial_\rho u = 0\,.\]
 This equation has the same form, up to an additional free parameter, as the advection equation in the original coordinates \eqref{eq:advection}, but the meaning of the coordinates is different. Solutions to the above equation in the bounded domain $\rho\in[0,1]$ correspond to solutions to the original advection equation in the unbounded domain $x \in[0,\infty)$. The free parameter $C$ expresses the freedom to prescribe the characteristic speed in the compactifying coordinate and the number of cycles in an infinite domain. To see this, we write the solution (\ref{eq:adv_sol}) in the new coordinates
%\[ u(\rho,\tau) = \sin\left(\frac{-\rho  \left(\tau ^2+2 \tau -1\right)+\tau ^2-1}{\sqrt{2 \rho ^2-2 \rho +1}+\rho  \tau +\rho -\tau } \right). \]
%\be\label{eq:adv_hypal} u(\rho,\tau) = -\sin\left(\frac{C(1-\rho)}{1+\rho} +\tau \right). \ee
\be\label{eq:hypsol} u(\rho,\tau) = -\sin \left(2\pi (C\,\Omega + \tau)\right)\,.\ee
 The solution is depicted in figure \ref{fig:red} at $\tau=0$ for two values of $C$. The number of cycles on the domain is finite and depends on $C$. The wave is resolved evenly through the compactified domain. In such representations of the solution it should be kept in mind that lines of constant $\tau$ do not correspond to lines of constant $t$. 

\begin{figure}[htb]
\center
\includegraphics[height=0.17\textheight]{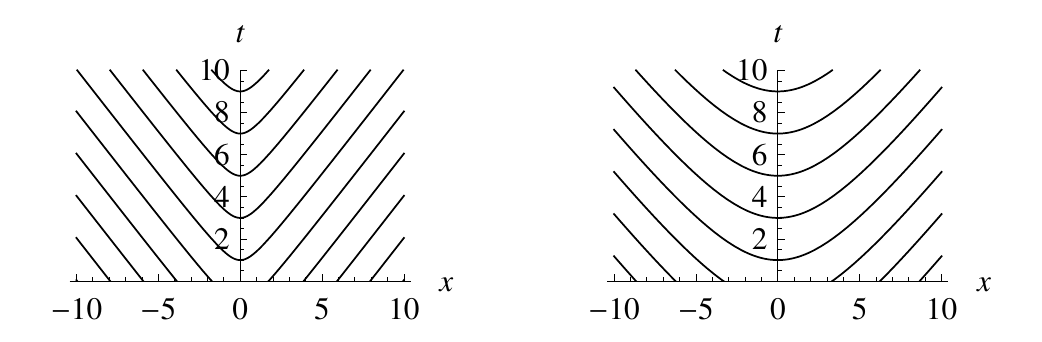}\hspace{1cm}
\caption{Rectangular hyperbola (\ref{eq:hyp1}) for $C=1$ and $C=5$ depicted on the $(x,t)$ plane. The hyperbolae are flatter near the origin for larger $C$ but their asymptotic behavior in leading order is independent of $C$.\label{fig:2}}
\end{figure}

The idea to introduce a coordinate transformation of time in combination with compactification comes from general relativity \cite{Penrose63}. The time function \eqref{tau} approaches characteristics of the advection equation asymptotically. In general relativity, infinity along characteristic directions is called \emph{null infinity}. Time functions whose level sets approach null infinity are called \emph{hyperboloidal} because their asymptotic behavior is similar to the asymptotic behavior of standard hyperboloids \cite{EardleySmarr,Friedrich83}. To see this, consider the rectangular hyperbola on the $(x,t)$ plane, $t^2-x^2 = C^2$, with a free parameter $C$. Shifting the hyperbola along the $t$ direction by $\tau$ gives $(t-\tau)^2 - x^2 =C^2$. Introducing $\tau$ as the new time coordinate we write
\be \label{eq:hyp1} \tau = t - \sqrt{C^2+x^2}. \ee
We plot in figure \ref{fig:2} two families of hyperbolae with $C=1$ and $C=5$. The asymptotes of these hyperbolae on $x>0$ are the characteristics of the advection equation, $\tau=t-x$, for any value of $C$. Time surfaces \eqref{tau} and \eqref{tauu} share the same asymptotic behavior, hence the name hyperboloidal for such surfaces. A suitable compactification along hyperboloidal surfaces sets the coordinate location of null infinity to a time independent value (see \cite{Zenginoglu:2007jw} for a discussion of conformal and causal properties of hyperboloidal surfaces and compactification on asymptotically flat spacetimes).

It is useful to summarize the hyperboloidal compactification technique with general expressions. We introduce coordinates $\tau$ and $\rho$ via
\be \label{eq:trafos} \tau = t - h(x), \qquad x = \frac{\rho}{\Omega}. \ee
The height function, $h$, must satisfy $|dh/dx| <1$ in units in which the speed of light (advection speed in the above example) is unity, so that $\tau$ is a time function. We also require that the gradient of the function $\Omega\equiv\Omega(\rho)$ is nonvanishing at its zero set. The zero set of $\Omega$ corresponds to infinity in $x$. The coordinate transformations have the Jacobian
\be\label{eq:jacob} \partial_t = \partial_\tau, \qquad \partial_x = -H \, \partial_\tau + \frac{\Omega^2}{L} \partial_\rho, \quad \textrm{where} \quad H:= \frac{dh}{dx}(\rho), \quad L:= \Omega - \rho \frac{d\Omega}{d\rho}. \ee
The advection equation becomes in this notation
\be \label{eq:tt} \partial_\tau u + \frac{\Omega^2}{(1-H)L} \partial_\rho u = 0. \ee
The specific choices \eqref{comp} and \eqref{tauu} satisfy
\be\label{condition}  \frac{\Omega^2}{(1-H)L} = \frac{1}{C}\,. \ee
In general, the time transformation must be chosen such that we have asymptotically in $x$, or equivalently as $\Omega$ approaches zero
\be\label{cond}1-H \sim O(\Omega^2).\ee
This relation formulates that the asymptotic boost of time surfaces needs to approach the speed of light as fast as the infinite compression of space to have a uniform outgoing characteristic speed in a compactified domain. This intuitive explanation suggests the names \emph{boost function} for $H$, and \emph{compress function} for $\Omega$.

\subsection{Wave equation}\label{sec:wave}
\begin{figure}[ht]
\center
\includegraphics[width=0.47\textwidth]{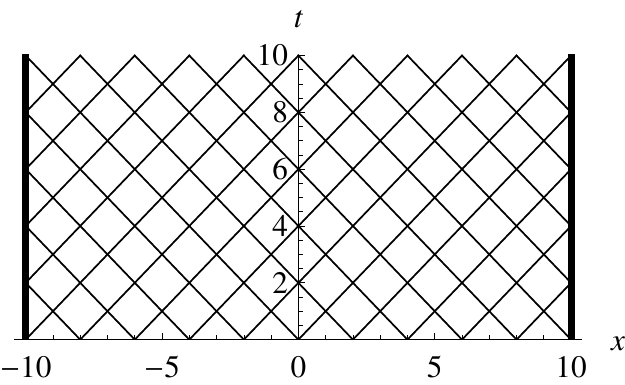}\hfill
\includegraphics[width=0.47\textwidth]{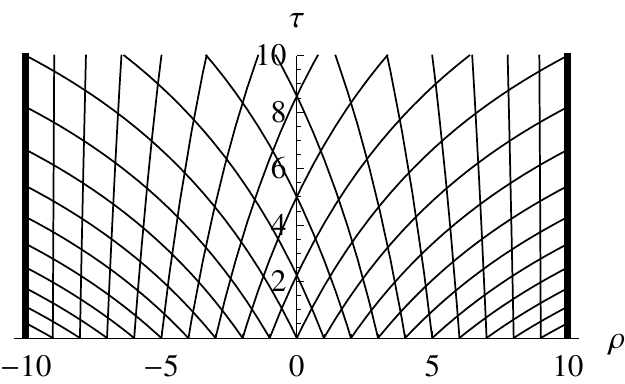}
\caption{On the left panel we plot the characteristics for the standard wave equation \eqref{eq:wave} on the bounded domain $x\in[-10,10]$. On the right panel hyperboloidal compactification has been applied with infinity located at $\rho=S =\pm10$. The standard wave equation on a bounded domain requires boundary conditions for the incoming characteristics from both boundaries. There are no incoming characteristics with hyperboloidal compactification.
\label{fig:wave_char} }
\end{figure}

The advection equation discussed in the previous section is a special example because its characteristics propagate in only one direction. A more representative example with incoming and outgoing characteristics is the wave equation,
\be \label{eq:wave} \partial_t^2 u -\partial_x^2 u=0. \ee
The characteristics of \eqref{eq:wave} on a bounded domain are plotted on the left panel of figure \ref{fig:wave_char}. We are interested in the problem on the unbounded domain $x \in (-\infty,\infty)$. With \eqref{eq:jacob} we get
\be\label{eq:wave_tr} \left[\partial_\tau^2 +\frac{\Omega^2}{(1-H^2) L} \left(2 H\partial_\tau\partial_\rho - \frac{\Omega^2}{L}\partial^2_\rho + \partial_\rho(H) \partial_\tau - \partial_\rho\left(\frac{\Omega^2}{L}\right)\partial_\rho \right) \right] u = 0. \ee
 The transformation \eqref{tauu} is not the right choice for the solution domain. Instead, we choose the height function of standard hyperboloids, $h(x) = \sqrt{S^2+x^2}$. We prescribe an arbitrary coordinate location for null infinity by modifying the zero set of the compress function. We set
%\be\label{eq:hyp2} \Omega = \frac{S^2-\rho^2}{2 \, C S} \quad \Rightarrow \quad L = \frac{S^2+\rho^2}{2 \, C S} \quad \mathrm{and} \quad H = \frac{2 S \rho}{S^2+\rho^2}. \ee
\be\label{eq:hyp2} \Omega =\frac{1}{2}\left(1-\frac{\rho^2}{S^2}\right) \quad \Rightarrow \quad L = \frac{1}{2}\left(1+\frac{\rho^2}{S^2}\right) \quad \mathrm{and} \quad H = \frac{2 S \rho}{S^2+\rho^2}. \ee

The parameter $S$ determines the coordinate location of null infinity on the numerical grid. The unbounded domain $x\in(-\infty,\infty)$ corresponds in the compactifying coordinate to $\rho\in(-S,S)$. The wave equation becomes
\be\label{eq:wave_comp} \partial_\tau^2 u + \frac{2\rho}{S}\partial_\tau\partial_\rho u -\Omega^2\,\partial_\rho^2 u + \frac{2 \,S\,\Omega}{S^2+\rho^2}\partial_\tau u + \frac{(3 S^2+\rho^2)\rho \,\Omega}{S^2 (S^2+\rho^2)} \partial_\rho u= 0. \ee
The equation evaluated at infinity, $\{\rho=\pm S\}$, takes the form 
\[\partial_\tau \left(\partial_\tau \pm 2\,\partial_\rho\right)\, u = 0, \]
reflecting that both boundaries are outflow boundaries and do not require boundary conditions.

The characteristics of the wave equation \eqref{eq:wave_comp} are plotted on the right panel of figure \ref{fig:wave_char}. No boundary conditions are needed because no characteristics enter the simulation domain. Furthermore, the outgoing characteristics leave the domain smoothly through the outflow boundaries. 

The initial value problem for the wave equation \eqref{eq:wave} is related to the initial value problem for \eqref{eq:wave_tr} by time evolution as indicated in figure \ref{fig:2}. This may be an undesired complication in practical applications. If we are interested in the evolution of compactly supported data, we can keep the $t$ surfaces in an interior domain which includes the data and apply hyperboloidal compactification only in the exterior domain as discussed in section \ref{sec:layer}.

The above calculation is made for the one dimensional, source free, constant coefficient wave equation. The method applies when variable coefficients or suitable lower order terms are present. The requirement on the variable coefficients in the principal part is that they are asymptotically constant. For the lower order terms we require a fall off behavior of at least $\Omega^2$ so that the division by $1-H^2$ in \eqref{eq:wave_tr} leads to a regular equation.
\subsection{Hyperbolic systems}\label{sec:sys}
Consider the linear, homogeneous system of partial differential equations with variable coefficients
\be \label{eq:sys} \partial_t \mathbf{u} = \mathbf{A} \partial_x \mathbf{u}, \ee
where $\mathbf{u} = (u_1(x,t), u_2(x,t), \dots, u_n(x,t))^T$, and $ \mathbf{A}$ is an $n\times n$ matrix that may depend on $x$. The transformation \eqref{eq:trafos} with Jacobian \eqref{eq:jacob} leads to
\be (\mathbbm{1} + H \mathbf{A}) \partial_\tau \mathbf{u} = \frac{\Omega^2}{L}  \mathbf{A} \partial_\rho \mathbf{u}. \ee
Assuming that the time transformation has been chosen to satisfy \eqref{cond}, we require that the polynomial remainder of $\det(\mathbbm{1} + H\mathbf{A})$ by $1-H$ vanishes asymptotically. This is a condition on the asymptotic form of the elements of $\mathbf{A}$. For example, taking $n=2$ we write
\[ \mathbf{A} = \left(\begin{array}{cc} a_{11} & a_{12} \\ a_{21} & a_{22} \end{array} \right)\,. \]
The asymptotic condition for the applicability of hyperboloidal compactification reads
\be \label{eq:cond} 1+ a_{11} + a_{22} - a_{12} a_{21} + a_{11} a_{22} = 0. \ee

A typical example is the wave equation (\ref{eq:wave}) written as a first order symmetric hyperbolic system. The wave equation takes the form (\ref{eq:sys}) in the auxiliary variables $v = \partial_t u$ and $w = \partial_x u$, with
\be \mathbf{u} = \left( \begin{array}{c} v \\ w \end{array} \right), \qquad \mathbf{A} = \left(\begin{array}{cc} 0 & 1 \\ 1 & 0 \end{array} \right). \ee
The condition (\ref{eq:cond}) is satisfied. The transformed system reads
\be\label{eq:wave_sys} \partial_\tau  \mathbf{u}  = \frac{\Omega^2}{(1-H^2)L}  \left(\begin{array}{cc} -H & 1 \\ 1 & -H \end{array} \right)\partial_\rho \mathbf{u}. \ee
Note that this equation is not the first order symmetric hyperbolic form of the transformed wave equation \eqref{eq:wave_tr}. The particular choice (\ref{eq:hyp2}) leads to the regular system
\[ \partial_\tau   \mathbf{u} = \frac{1}{2\,S^2}  \left(\begin{array}{cc} -2 S \rho & S^2+\rho^2 \\ S^2+\rho^2 & -2 S \rho \end{array} \right)\partial_\rho \mathbf{u}. \]
The outer boundaries at $\rho=\pm S$ are pure outflow boundaries.

As a further example consider the one dimensional Maxwell equations for the electric and magnetic fields $(\bar{E},\bar{H})$
\[ \partial_t \bar{E} = - \frac{1}{\epsilon}\partial_x \bar{H}, \quad  \partial_t \bar{H} = - \frac{1}{\mu}\partial_x \bar{E}, \]
The electric permittivity $\epsilon$ and the magnetic permeability $\mu$ may be point-dependent. The equations have the form (\ref{eq:sys}) with
\[ \mathbf{u} = \left( \begin{array}{c} \bar{E} \\ \bar{H} \end{array} \right), \quad \textrm{and} \quad \mathbf{A} = -\frac{1}{\epsilon\mu}\left(\begin{array}{cc} 0 & \mu \\ \epsilon & 0 \end{array} \right). \]
We get after hyperboloidal compactification
\be\label{eq:maxw} \partial_\tau  \mathbf{u}  = -\frac{\Omega^2}{(\epsilon\mu-H^2)L}  \left(\begin{array}{cc}  H & \mu \\ \epsilon & H \end{array} \right)\partial_\rho \mathbf{u}. \ee
In vacuum outside a compact domain we have $\epsilon=\epsilon_0$ and $\mu=\mu_0$, where $\epsilon_0$ and $\mu_0$ are the electric and the magnetic constants. We need to choose the asymptotic behavior of $H$ such that $\sqrt{\epsilon_0 \mu_0} - H \sim O(\Omega^2)$. Then the Maxwell equations behave similarly to the wave equation (\ref{eq:wave_sys}) near null infinity.

The example of Maxwell equations suggests that including lower order terms or variable characteristic speeds in a compact domain are straightforward in the hyperboloidal method as long as the asymptotic form of the equations are suitable. The asymptotic characteristic speeds need to be constant and lower order terms need to have compact support or fall off sufficiently fast. Hyperboloidal compactification can then be applied outside a compact domain as discussed in the next section.

\subsection{Hyperboloidal layers}\label{sec:layer}
It may be desirable to employ specific coordinates in a compact domain without the time transformation or the compactification required by the hyperboloidal method. One reason is technical. Elaborate numerical techniques to deal with shocks, scatterers, and media assume predominantly specific coordinates. It may be impractical to modify these methods to work with hyperboloidal compactification throughout the simulation domain. Another reason is initial data. One may be interested in the evolution of certain (compactly supported) initial data prescribed on a level set of $t$. Therefore, it may be favorable to restrict the hyperboloidal compactification to a layer. 

We discuss briefly the perfectly matched layer (PML)  by B\'erenger \cite{Berenger} to set the stage for hyperboloidal layers.  In the PML method one attaches an absorbing medium---a layer--- to the domain of interest such that the interface between the interior domain and the exterior medium is transparent independent from the frequency and the angle of incidence of the outgoing wave. Inside the layer the solution decays exponentially in the direction perpendicular to the interface. As a consequence, the solution is close to zero at the outer boundary of the layer where any stable boundary condition may be applied. The reflections from the outer boundary may be ignored if the layer is sufficiently wide.

The success of the PML method lies in the transparency of the interface between the interior domain and the layer. This property finds explanation in the interpretation of Chew and Weedon of the PML as the analytic continuation of the equations into complex coordinates \cite{ChewWeed}. The challenge is then to find suitable choices of the equations and the free parameters that lead to exponential damping of the solution in a stable way, which may be difficult depending on the problem \cite{AbarGott, AppHagKr, Becache}. 

\begin{figure}[ht]
\center
\includegraphics[height=0.17\textheight]{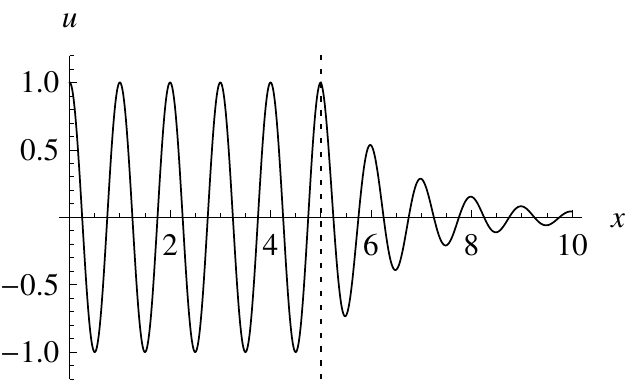}
\caption{The plane wave solution with the spatial coordinate analytically extended into the complex plane \eqref{eq:pml} demonstrates how outgoing waves are damped exponentially in the absorbing layer of the PML method. We set $\sigma=0.1$ and $R=5$. The dashed line indicates the location of the interface.
\label{fig:pml}}
\end{figure}

A simple example omitting details of the method beyond our needs demonstrates the basic idea. We perform an analytic continuation of $x$ into the complex plane beyond a certain interface $R$. The coordinate $x$ is written in terms of its real and imaginary parts as $\mathrm{Re}(x) + i \sigma \mathrm{Im}(x) \Theta(x-R)$, where $\sigma$ is a positive parameter, and $\Theta$ denotes the Heaviside step function. Setting $k=1$, the plane wave \eqref{plane} at time $t=0$ becomes
\be\label{eq:pml} u(x,0)=e^{ 2 \pi i \mathrm{Re}(x)} e^{-2 \pi \sigma \mathrm{Im}(x) \Theta(x-R)}\,. \ee
The strength of the exponential decay is controlled by a free parameter $\sigma$. The solution is plotted in figure \ref{fig:pml} for $\sigma=0.1$. 

For a hyperboloidal layer we perform a real coordinate transformation, both of space and time, beyond a certain timelike surface $x=R$ called the interface. We set
\be\label{eq:hl1x} x - R = \frac{\rho-R}{\Omega}, \qquad \Omega = 1-\frac{(\rho-R)^2}{(S-R)^2}\Theta(\rho-R)\,,\ee
where the coordinate location of infinity satisfies $S>R$. The width of the layer is $S-R$. This transformation is partly motivated by the agreement between the coordinates $x$ and $\rho$ along the interface in the sense that 
\be\label{eq:agree} x(R)=R, \quad \frac{d x}{d\rho} (R)= 1, \quad \quad \frac{d^2 x}{d\rho^2} (R)= 0\,.\ee

\begin{figure}[ht]
\center
\includegraphics[width=0.91\textwidth]{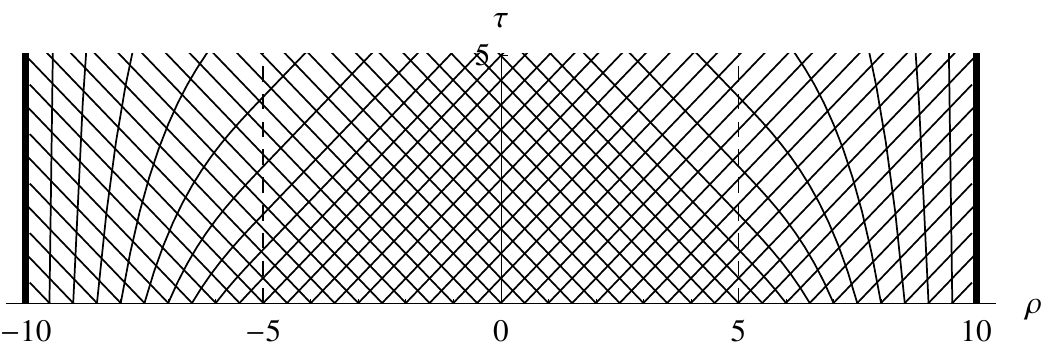}
\caption{Characteristic structure for hyperboloidal layers with boundaries located at $\rho=\pm S=\pm 10$ and interfaces at $\rho=\pm R=\pm 5$. Compare the inner domain $[-5,5]$ to the left panel of figure \ref{fig:wave_char}, and the layers $[-10,-5]$ and $[5,10]$ to the right panel of figure \ref{fig:wave_char}. No characteristics enter the computational domain.
\label{fig:layer_chars}}
\end{figure}

A simple prescription for the boost function is to require unit outgoing characteristic speed across the layer. For example, the outgoing and incoming characteristic speeds $c_\pm$ for the wave equation on $\rho>0$ are
\be\label{eq:chars} c_{\pm} = \frac{\Omega^2}{L (\pm 1 - H)}, \qquad \mathrm{with} \qquad L = \Omega - (\rho-R)\frac{d\Omega}{d\rho}.  \ee
The requirement of unit outgoing characteristic speed reads $c_+=1$, implying
\be\label{eq:hl1om} 1-H = \frac{\Omega^2}{L}. \ee
For $\rho<-R$ we require $c_-=-1$. The resulting characteristics are depicted in figure \ref{fig:layer_chars}.  For $|\rho|<R$ we obtain the standard characteristics in $(x,t)$ coordinates. For $|\rho|>R$ we obtain the hyperboloidal characteristics in $(\rho,\tau)$ coordinates (compare figure \ref{fig:wave_char}). %No characteristics enter the computational domain.

\begin{figure}[ht]
\center
\includegraphics[height=0.17\textheight]{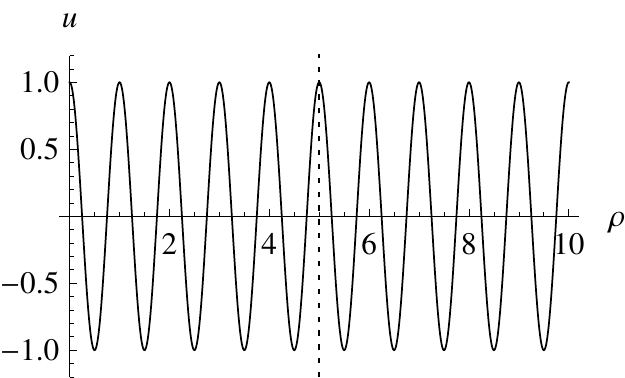}\hspace{1cm}
\includegraphics[height=0.17\textheight]{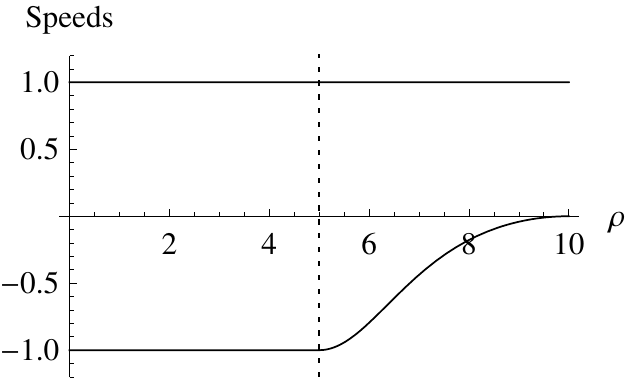}
\caption{On the left panel we plot the plane wave solution with the compress function given in \eqref{eq:hl1x} and the boost function given in \eqref{eq:hl1om}. The characteristic speeds are depicted on the right panel in which the top and the bottom curves correspond to outgoing and incoming characteristic speeds. The domain is given by $\rho\in[0,10]$ with a compactifying layer starting at $\rho=R=5$ as indicated by the dashed lines. The incoming characteristic speed at the outer boundary vanishes, therefore no outer boundary conditions are needed in the layer.
\label{fig:hl1}}
\end{figure}

The plane wave in the layer has the same form as in the interior by \eqref{eq:hl1x} and \eqref{eq:hl1om} (figure \ref{fig:hl1}; the dashed line indicates  the interface). The important difference to the standard method is that the incoming characteristic speed vanishes at the outer grid boundary.

\begin{figure}[ht]
\center
\includegraphics[height=0.17\textheight]{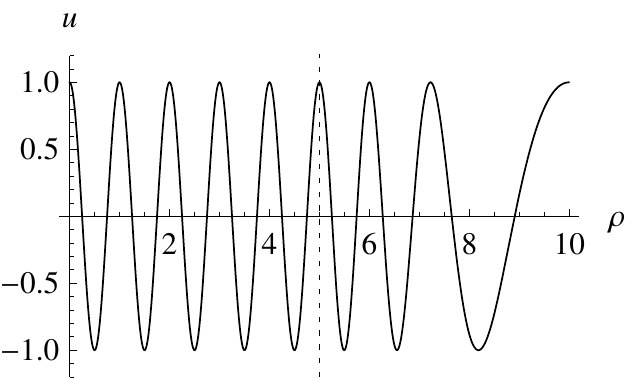}\hspace{1cm}
\includegraphics[height=0.17\textheight]{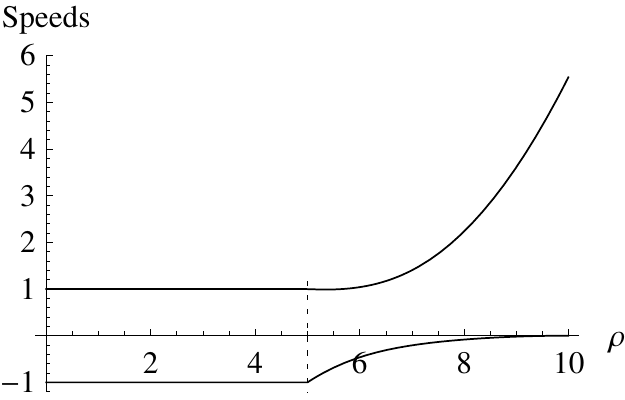}
\caption{On the left we plot the plane wave solution with \eqref{eq:hl2om} and \eqref{eq:hlh} for $C=3$ in the otherwise same setting as in figure \ref{fig:hl1}. There are fewer oscillations in the layer due to the stronger redshift causing a higher outgoing characteristic speed as depicted on the right panel.
\label{fig:hl2}}
\end{figure}

There is a large freedom in the choice of compress and boost functions, which may be exploited for specific purposes. As an example, take a linear compress function
\be\label{eq:hl2om} \Omega = 1- \frac{\rho-R}{S-R} \Theta(\rho-R), \ee
and the boost function of standard hyperboloids translated by $R$
\be\label{eq:hlh} H = \frac{x-R}{\sqrt{(x-R)^2+C^2}}\Theta(x-R)\, . \ee
The resulting representation of the plane wave is plotted on the left panel of figure \ref{fig:hl2} for $C=3$. Now we have fewer oscillations in the layer than in the interior because of the spatial redshift controlled by $C$. A strong redshift, and consequently few spatial oscillations, may be preferable in a high order space discretization scheme in which a few grid points per wavelength are sufficient for good accuracy. On the right panel of figure \ref{fig:hl2}, however, we see that strong spatial redshift comes at a price: the outgoing coordinate speed increases strongly in the layer. Evaluation of $c_+$ from \eqref{eq:chars} gives at infinity $c_+ = 2 (S-R)^2/C^2$. A small value for $C$ with a wide layer leads to a high outgoing characteristic speed, which requires small time steps in an explicit time integration algorithm. 

The balance between accuracy in time and in space is influenced by the compress and boost functions, and the free parameters included in them. The specific choices will depend on the requirements of the problem. It can be expected that in most applications the outgoing characteristic speed in the layer will be chosen close to the speed in the interior domain as in figure \ref{fig:hl1}.

The hyperboloidal layer is similar to the PML in that both methods allow us to solve the  equations of interest in standard coordinates in an inner domain. In both methods the interface between the inner domain and the layer is transparent on the analytical level, independent of the frequency and angle of incidence of the outgoing wave. The key difference is that the PML absorbs the outgoing wave so that it is damped exponentially, whereas the compactifying layer transports it to infinity. The solution in the hyperboloidal layer is of interest, as opposed to the solution in the absorbing layer. In fact, the solution at the outer boundary of the layer is of special interest in radiative systems because it gives the radiation signal as measured by an idealized observer at infinity. 

\section{Multiple dimensions}\label{sec:mult}
Hyperboloidal compactification is applicable in multiple dimensions when compactification is performed in the outgoing direction. The main difference in multiple dimensions is that a rescaling of the unknown fields needs to be performed such that the fields attain a non-vanishing finite limit at infinity. The rescaling depends on the fall off behavior of these fields, and therefore on the equation and the spatial dimension. 
\pagebreak

\subsection{Rescaling}\label{sec:resc}
We begin our discussion with the three dimensional wave equation on flat spacetime
\[ \left(-\partial_t^2 + \Delta_{\mathbb{R}^3}\right) u = 0. \]
Here, $\Delta_{\mathbb{R}^3}$ is the Laplace operator on three dimensional Euclidean space. We write the wave equation in spherical coordinates to single out the outgoing direction
\[ \left(-\partial_t^2 + \partial_r^2 + \frac{2}{r} \partial_r + \frac{1}{r^2} \Delta_{\mathbb{S}^2} \right) u = 0,   \]
with $\Delta_{\mathbb{S}^2}$, the Laplace--Beltrami operator on the two sphere, and $r$, the radial coordinate. In section \ref{sec:wave} we showed that hyperboloidal compactification introduces a divisor that is proportional to the square of the compress function. With the compactification $r=\rho/\Omega$ we see that any term beyond the flat wave operator on the $(r,t)$ plane that falls off as $r^{-2}$ or faster is multiplied with at least $\Omega^2$, and is therefore regular under hyperboloidal compactification. 

The angular part of the wave equation in spherical coordinates admits a regular compactification. The first radial derivative term, however, results in a singular operator at infinity. We resolve this problem by rescaling. The three dimensional wave equation on the $(r,t)$ plane takes the form of the one dimensional wave equation in the rescaled variable $v:=r u$. We get \[\left(-\partial_t^2 + \partial_r^2 + \frac{1}{r^2} \, \Delta_{\mathbb{S}^2} \right) v = 0,\]
which admits a regular hyperboloidal compactification.

This procedure generalizes to dimensions other than three. The essential feature of the rescaling is that it takes care of the fall off behavior of the scalar field. In three dimensions $r u$ attains a regular limit at infinity; in $n$ dimensions the rescaling depends on $n$. Solutions to the wave equation decay asymptotically as $r^{-(n-1)/2}$ by energy conservation. We expect therefore that the $n$ dimensional wave equation written for the rescaled variable $v := r^{(n-1)/2} u$ admits a regular hyperboloidal compactification.

The $n$ dimensional wave equation in spherical coordinates reads 
\[  \left(-\partial_t^2 +  \partial_r^2 + \frac{n-1}{r}\partial_r + \frac{1}{r^2} \Delta_{\mathbb{S}^{n-1}} \right) u = 0, \] 
where $\Delta_{\mathbb{S}^{n-1}}$ is the Laplace-Beltrami operator on the $n-1$ sphere. Transforming to the rescaled variable $v := r^{(n-1)/2} u$, we get
\[  \left(-\partial_t^2 +  \partial_r^2  - \frac{1}{4r^2} (n-1)(n-3) +  \frac{1}{r^2} \Delta_{\mathbb{S}^{n-1}}  \right) v = 0. \]
All terms beyond the one dimensional wave operator fall off as $r^{-2}$ and are therefore amenable to a regular hyperboloidal compactification. 

The inclusion of sources or of suitable nonlinearities is straightforward under the condition that the corresponding terms in the equation fall off sufficiently fast. For example, a power nonlinearity on the right hand side of the wave equation of the type $u^p$ leads to the source term $v^p r^{-(n-1)(p-1)/2}$. This term is regular under hyperboloidal compactification if the power of $r$ is $-2$ or lower, which implies $p\geq1+4/(n-1)$. The critical power for which equality is satisfied is also the critical conformal power for semilinear wave equations with a power nonlinearity. This apparent coincidence is explained in the next section. 

\subsection{Conformal method}\label{sec:conf}
Compactification of spacetimes with a suitable time transformation as proposed by Penrose in \cite{Penrose63} as well the hyperboloidal initial value problem as proposed by Friedrich \cite{Friedrich83} employ conformal methods. The conformal language is prevalent in studies of spacetimes in general relativity \cite{Wald84}. In this section we discuss the hyperboloidal compactification using  conformal techniques. This viewpoint is of theoretical and practical interest because it reveals the interplay of conformal geometry, partial differential equations, and numerical methods within the hyperboloidal approach, and also simplifies the implementation of the method for certain problems. 

In general, a hyperboloidal time transformation with a spatial compactification leads to a singular metric. Consider the Minkowski metric in spherical coordinates
\be\label{eq:mink} \eta = - dt^2 +r^2 dr^2 + r^2d\sigma^2,   \ee
where $d\sigma^2$ is the standard metric on the unit sphere. Introducing new coordinates $\tau$ and $\rho$ as in \eqref{eq:trafos} gives \cite{Zenginoglu:2007jw}
\[ \eta = -d\tau^2 - \frac{2 H L}{\Omega^2} d\tau d\rho + \frac{1-H^2}{\Omega^4} L^2 d\rho^2 + \frac{\rho^2}{\Omega^2} d\sigma^2. \]
This representation of the Minkowski metric is singular at infinity. The singularity is removed by conformally rescaling the metric, 
\be\label{conf_met} g = \Omega^2 \eta = - \Omega^2\,d\tau^2 - 2 H L\, d\tau d\rho + \frac{1-H^2}{\Omega^2} L^2 d\rho^2 + \rho^2\, d\sigma^2. \ee
The conformal metric $g$ is regular at null infinity by \eqref{cond}, and can be extended beyond null infinity in a process referred to as conformal extension \cite{Penrose63,Wald84, Penrose65}. In this context, the function $\Omega$ is called the conformal factor. The zero set of the conformal factor corresponds to null infinity where it has a non-vanishing gradient. These properties of the conformal factor underlie our choices for the compress function in the previous sections.

Partial differential equations within the conformal framework have first been studied for fields with vanishing rest mass, such as scalar, electromagnetic, and gravitational fields \cite{Penrose65}. The key observation is that the complicated asymptotic analysis of solutions to these partial differential equations is replaced with local differential geometry by considering the conformally transformed equations in a conformally extended, regular spacetime \cite{Geroch77, Friedrich02}. 

We discuss the wave equation on Minkowski spacetime as an example. Under a conformal rescaling of the Minkowski metric, $g=\Omega^2\,\eta$, the wave equation transforms as \cite{Wald84, Penrose65}
\be \label{conf_invar} \left(\Box_g - \frac{n-1}{4n}R[g]\right)
v = \Omega^{-(n+3)/2} \, \Box_\eta\, u,\quad\textrm{with} \quad v :=
\Omega^{(1-n)/2} u.\ee Here, $\Box_g :=
g^{\mu\nu}\nabla_\mu\nabla_\nu$ is the d'Alembert operator with
respect to $g$, $R[g]$ is the Ricci scalar of $g$, and $n$ is the
spatial dimension of the spacetime. The rescaling with $\Omega$ in the definition of the variable $v$ is asymptotically equivalent to the rescaling in Section \ref{sec:resc}, where we factor out the fall off behavior of $u$ such that the rescaled variable $v$ has a non-vanishing limit at null infinity. To see this, consider a specific choice for $\Omega$ from previous sections, say $\Omega(\rho)=1-\rho$ as in \eqref{comp}. In terms of the coordinate $r=\rho/\Omega$, the conformal factor reads $\Omega(r) = (1+r)^{-1}$, which behaves  asymptotically as $r^{-1}$. Therefore the definition of $v$ in \eqref{conf_invar} corresponds asymptotically to the definition of $v$ in section \ref{sec:resc}.

We can also explain the observation made at the end of section \ref{sec:resc} concerning the agreement between the critical conformal power and the critical power for which hyperboloidal compactification leads to a regular equation. Using $\Box_\eta u = u^p$ and the definition of $v$ in \eqref{conf_invar} we get at the right hand side of the conformally invariant wave equation $\Omega^{\left((n-1) p - (n+3)\right)/2}v^p$. Regularity of this term at infinity, where the conformal factor vanishes, requires $(n-1) p \geq n+3$, which is the condition of Section \ref{sec:resc}. Equality is obtained for $p_c = 1+4/(n-1)$ for which the semilinear wave equation is conformally invariant, hence $p_c$ is called the critical conformal power.

The conformal approach may be useful for various reasons. It extends directly to asymptotically flat backgrounds with non-vanishing curvature \cite{Zenginoglu:2007jw, Zenginoglu:2008wc}. It also helps identifying the transformation behavior of the equations independent of coordinates. For example, Yang-Mills and Maxwell equations are conformally invariant, and therefore do not require a rescaling of the variables (the specific variables in which the covariant equations are written may not be conformally invariant and may require a rescaling).

Hyperboloidal compactification with rescaling, as presented in section \ref{sec:resc}, seems numerically feasible in any spatial dimension. In even spatial dimensions, however, the notion of conformal infinity may not be feasible \cite{HollWald}. It is an open question whether the hyperboloidal technique applies to this case. If the method fails in even dimensions, it should be interesting to understand whether this failure is related to the violation of Huygens' principle.

\subsection{Hyperboloidal layers in multiple dimensions}\label{sec:mult_lay}
In multiple dimensional problems, the layer technique is employed along the outgoing direction $r$. We employ a compactifying coordinate $\rho$ defined via $r=\rho/\Omega$. The compress function needs to be unity in a compact domain bounded by radius $R$, be sufficiently smooth across the interface at $r=\rho=R$, and vanish at a coordinate location $S>R$ with non-vanishing gradient. A suitable choice is
\[ \Omega = 1- \left(\frac{\rho-R}{S-R}\right)^4 \Theta(\rho-R), \qquad  L = 1 + \frac{(\rho-R)^3 (3\rho+R)}{(S-R)^4} \Theta(\rho-R)\,.\]
The coordinates $r$ and $\rho$ coincide up to second order along the interface  (compare \eqref{eq:agree}). The height function in the new time coordinate depends only on the radial coordinate: $\tau = t-h(r)$. The boost function is given by $H=dh/dr$. It may be set such that the outgoing characteristic speed through the layer is unity. For the wave equation this corresponds to
\[ 1-H = \frac{\Omega^2}{L}, \]
as in \eqref{eq:hl1om}. Finally, the rescaling of the variable is performed with the compress function noting that asymptotically $\Omega \sim r^{-1}$. Further details are presented in section \ref{sec:three} in which the layer technique is applied to solve the wave equation numerically in three spatial dimensions.

The extension of the hyperboloidal method from one dimension to multiple dimensions is fairly straightforward in spherical coordinates, which may be, however, too restrictive for the grid geometry in applications. The discussion of the conformal regularity of the Minkowski metric suggests that the method can also be applied in nonspherical coordinates as discussed in the next section.

\subsection{Nonspherical coordinate systems}\label{sec:nonsph}
Spherical coordinates are suitable for hyperboloidal compactification because cuts of null infinity have spherical topology \cite{Penrose65}. The coordinate shape of null infinity, however, does not have to be a sphere. We can employ any coordinate system with a unique outgoing direction for the compactification. 

An example for a nonspherical coordinate system, useful especially in electromagnetism, is the prolate spheroidal coordinate system. The relation between Cartesian coordinates $\{x,y,z\}$ and prolate spheroidal coordinates $\{\mu,\nu,\phi\}$ reads
\[ x=\sinh \mu \sin \nu \cos \varphi,\qquad y = \sinh\mu \sin\nu \sin\varphi,\qquad z = \cosh \mu \cos\nu. \]
We have $r^2=\sinh^2\mu+\cos^2\nu$. Here, $\mu$ is the outgoing direction that has closed coordinate surfaces. Compactification is performed along $\mu$. Using the conformal method, we argue that if conformal compactification of Minkowski spacetime in these coordinates leads to a regular metric, the equations we solve on that background will be regular. The Minkowski metric reads
\[ \eta = -dt^2+dx^2+dy^2+dz^2 = -dt^2+(\sinh^2\mu + \sin^2\nu)\, 
(d\mu^2+d\nu^2) +\sinh^2\mu\sin^2\nu\,d\varphi^2.\]
We introduce a new time coordinate by setting
\[ \tau = t - \sqrt{1+\sinh^2\mu}.\]
The metric becomes
\[ \eta = -d\tau^2-2\sinh\mu\,d\mu d\tau+\sin^2\nu\,d\mu^2+(\sinh^2\mu + \sin^2\nu)d\nu^2+ \sinh^2\mu \sin^2\nu d\varphi^2. \]
Compactification along $\mu$ is performed via
\[ \sinh\mu = \frac{2\rho}{1-\rho^2} = \frac{\rho}{\Omega}, \qquad  d\mu = \frac{d\rho}{\Omega}.\]
The conformal metric $g=\Omega^2\eta$ becomes
\[ g = -\Omega^2d\tau^2-2\rho\, d\rho \,d\tau+\sin^2\nu \,d\rho^2+(\rho^2+\Omega^2\sin^2\nu)d\nu^2+\rho^2\sin^2\nu\,d\varphi^2. \]
The qualitative behavior of this metric near infinity is similar to the conformal Minkowski metric given in \eqref{conf_met}; the only difference is the coordinate representation. Therefore we conclude that hyperboloidal compactification of suitable hyperbolic equations in prolate spheroidal coordinates leads to regular equations as for spherical or oblate spheroidal coordinates. 

It is an open question whether the hyperboloidal method applies to Cartesian or cylindrical coordinates. The difficulty is to ensure regularity of the equations at corners and along edges with respect to limits to infinity. This requirement implies a restriction on the geometry of both the interface boundary and the numerical outer boundary. Even with the requirement of smooth coordinate surfaces, it would be useful to extend the method such that the layer has an arbitrary shape in coordinate space. This generalization would increase the efficiency of the method when dealing with scatterers with irregular shape. 

\section{Numerical experiments}\label{sec:num}
The aim of this chapter is to demonstrate numerical implementations of hyperboloidal compactification on simple examples in one and three spatial dimensions. The application of the method in two dimensions is an open question as discussed in \ref{sec:conf} and is left for future work. Also left for future work is an in-depth numerical analysis of the method in comparison with other boundary treatments, such as absorbing boundary conditions or perfectly matched layers.

In previous sections, we have seen that we may employ an arbitrary coordinate system in an interior domain restricting compactification to a layer outside that domain. This idea is based on the matching method with a transition zone presented in \cite{Zenginoglu:2007jw}. Numerical applications of the matching, however, are inefficient due to a blueshift in frequency in the matching region and many arbitrary parameters for the transition function \cite{Zenginoglu:2009ey, Zenginoglu:2010zm}. A hyperboloidal layer gives a sufficiently smooth interface without a transition function. 

In this section, we compare the numerical accuracy of solutions with and without the layer in one dimension. In three dimensions we focus on the application of hyperboloidal layers for the wave equation. Calculations using hyperboloid foliations, that is, constant mean curvature surfaces without the layer, have been presented in \cite{Zenginoglu:2010zm}.

\subsection{One spatial dimension}\label{sec:one}
\subsubsection{Analytical setup}
Consider the Maxwell equations \eqref{eq:maxw}. Assume that the electric permittivity and the magnetic permeability are constant and have unit value. The transformed Maxwell system for the unknown vector $\mathbf{u} = (\bar{E}, \bar{H})^T$ reads
\be\label{num_maxw} \partial_\tau  \mathbf{u}  = -\frac{\Omega^2}{(1-H^2)L}  \left(\begin{array}{cc}  H & 1 \\ 1 & H \end{array} \right)\partial_\rho \mathbf{u}. \ee
The characteristic speeds are $c_\pm=-\Omega^2/((\pm 1+H)L)$. This system is similar to the wave equation written in first order symmetric hyperbolic form \eqref{eq:wave_sys}, so our results apply both to Maxwell and wave equations in one dimension. 

We experiment with two sets of the compress and boost functions. First we employ the hyperboloid foliation everywhere in the simulation domain. We set as in \eqref{eq:hyp2}
\[ \Omega = \frac{1}{2} \left(1-\frac{\rho^2}{S^2}\right) \quad \mathrm{and} \quad H = \frac{2 S \rho}{S^2+\rho^2}. \]
The characteristic speeds read $c_\pm=\pm(1\pm\rho/S)^2/2$. The minus sign corresponds to the incoming speed at the right boundary, which vanishes at $\rho=S$. The plus sign corresponds to the incoming speed at the left boundary, which vanishes at $\rho=-S$. The qualitative behavior of the characteristics is the same as depicted on the right panel of figure \ref{fig:wave_char}. The time step in an explicit time integration scheme is restricted by the maximum absolute value of the characteristic speed which reads $|c|_{\rm max} = 2$.

The second set of compress and boost functions are for hyperboloidal layers. We set the compress function as in \eqref{eq:hl1x}
\[ \Omega = 1- \frac{(\rho-R)^2}{(S-R)^2} \Theta(\rho-R). \]
We determine the boost function from the requirement of unit outgoing characteristic speeds through the layers. We set
\[ H = 1-\frac{\Omega^2}{L} \quad \mathrm{for} \quad \rho>R\,, \qquad \mathrm{and} \quad H = -1 - \frac{\Omega^2}{L}   \quad \mathrm{for} \quad \rho<-R\,.\]
The characteristic structure for the resulting equations are depicted in figure \ref{fig:layer_chars}.

\subsubsection{Numerical setup}
The hyperboloidal method is essentially independent of numerical details due to its geometric origin. We discretize \eqref{num_maxw} employing common methods: an explicit 4th order Runge--Kutta time integrator and finite differencing in space with 4th, 6th, and 8th order accurate centered operators. At the boundaries we apply one-sided stencils of the same order as the inner operator. We use Kreiss--Oliger type artificial dissipation to suppress numerical high-frequency waves \cite{KreissOliger73}. For a $2p-2$ accurate scheme we choose a dissipation operator $D_{\mathrm{diss}}$ of order $2p$ as
\[ D_{\mathrm{diss}} =\ \epsilon\, (-1)^p \frac{h^{2p-1}}{2^p}\,D_+^p D_-^p,  \]
where $h$ is the grid size, $D_\pm$ are forward and backward
finite differencing operators and $\epsilon$ is the dissipation
parameter. 

Both for the hyperboloid foliation and the layer we set $S=10$. The simulation domain is then given by $\rho\in[-10,10]$, which corresponds to the unbounded domain $x\in(-\infty,\infty)$. The interface for the layer is at $R=\pm5$. Hyperboloidal coordinates $(\rho,\tau)$ coincide with standard coordinates $(x,t)$ within the domain $|\rho|\leq 5$.

We solve the initial value problem for \eqref{num_maxw} with a Gaussian wave packet centered at the origin for the electric field and vanishing data for the magnetic field. We set at the initial time surface
\[ \bar{E}(\rho,0) = e^{-\rho^2}, \quad \bar{H}(\rho,0) = 0\,.\]
The solutions with respect to the hyperboloid foliation and the hyperboloidal layer correspond to different initial value problems due to the different time surfaces. The solution constructed using the hyperboloidal layer corresponds, however, to the solution that we would obtain using the standard coordinates $(x,t)$.

\subsubsection{Results}
\begin{figure}[ht]
\center
\includegraphics[height=0.19\textheight]{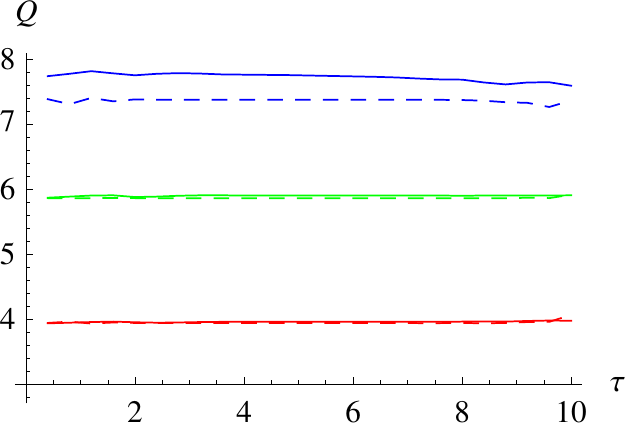}
\caption{Convergence factors in time in the $L_2$ norm for the Maxwell equations in one spatial dimension. Solid curves represent the hyperboloid foliation, dashed curves the layer method. The curves indicate from bottom to top 4th, 6th, and 8th order convergence in accordance with the order of implemented finite differencing operators.  
\label{fig:conv1d}}
\end{figure}

Figure \ref{fig:conv1d} shows convergence factors for the electric field from a three level convergence test with 100, 200, and 400 grid cells. The convergence factor $Q$ is calculated for the electric field in the $L_2$ norm via $Q:=\log_2\frac{ \|\bar{E}^{low}-\bar{E}^{med}\|}{\| \bar{E}^{med}-\bar{E}^{high}\|}$. The factors are in accordance with the implemented finite difference operators.

\begin{figure}[ht]
\center
\includegraphics[height=0.19\textheight]{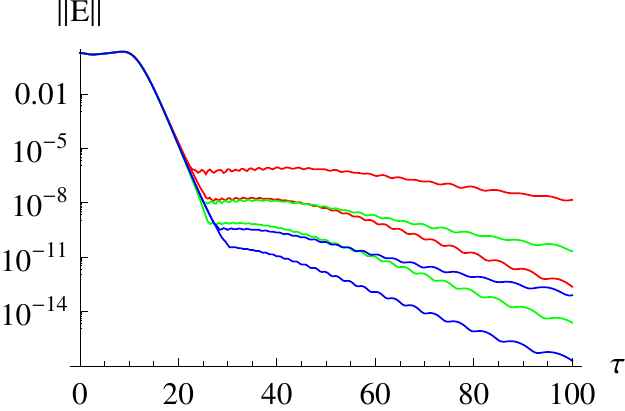}\hspace{7mm}
\includegraphics[height=0.19\textheight]{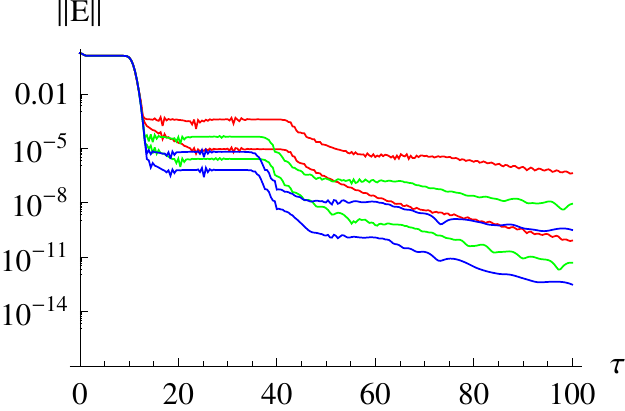}
\caption{The $L_2$ norm of the electric field in time for 4th (red), 6th (green), and 8th (blue) order finite differencing with respect to the hyperboloid foliation (left) and the hyperboloidal layer (right). The curve indicating smaller errors to each order corresponds to the solution with artificial dissipation. The plot shows that the effect of dissipation to the quality of the solution is comparable to using a higher order numerical discretization. Errors with the layer method are larger than errors with the hyperboloid foliation.
\label{fig:refl}}
\end{figure}

The total energy is radiated to infinity leaving the zero solution behind in accordance with Huygens' principle for Maxwell equations in one dimensional flat spacetime. Therefore, a good measure of the quality of the boundary treatment is the value of the unknown after the initial wave packet leaves the simulation domain. This value is related to the numerical reflection coefficient of the boundary. The analytical reflection coefficients at the interfaces and at the boundaries are zero in the hyperboloidal method. Numerically, however, the interfaces and the finite differencing at the outer boundaries cause reflections. Consider the $L_2$ norm of the solution as a function of time plotted in figure \ref{fig:refl}. The plots depict the field with respect to the hyperboloid foliation (left) and with respect to the layer (right). The errors are larger for the layer. Red, green, and blue curves represent solutions calculated with 4th, 6th, and 8th accurate finite difference operators, respectively. To each order the solution is calculated with and without dissipation. We observe that artificial dissipation strongly reduces numerical errors. Figure \ref{fig:refl} shows that the error at late times with 4th order stencils with dissipation has the same order of magnitude as the one with 8th order stencils without dissipation.

\begin{figure}[ht]
\center
\includegraphics[height=0.19\textheight]{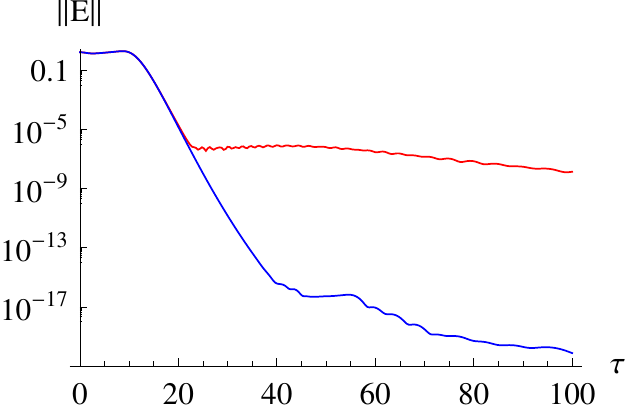}\hspace{7mm}
\includegraphics[height=0.19\textheight]{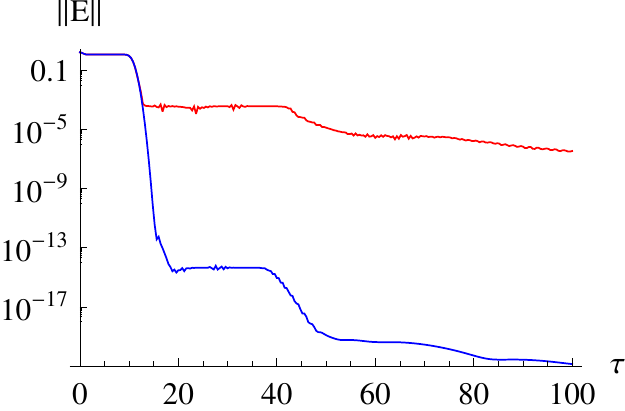}
\caption{Comparison of numerical errors as indicated by late time values of the $L_2$ norm of the electric field between a numerical solution with low (red) and high accuracy (blue) for the hyperboloid (left) and the layer (right) method. Accuracy with hyperboloidal compactification is not restricted by the boundary treatment.
\label{fig:compare}}
\end{figure}

The main result from figures \ref{fig:conv1d} and \ref{fig:refl} is that the numerical error is reduced drastically by higher resolution, higher order discretization, and artificial dissipation---without changing the boundary treatment. The boundary treatment does not introduce errors into the solution that are put in by hand. To emphasize this point, we plot in figure \ref{fig:compare} two solutions using the hyperboloid foliation and the hyperboloidal layer. For each method, one solution is obtained with 4th order finite differencing, without dissipation, and 100 grid cells (red curve); the other solution is obtained with 8th order finite differencing, with dissipation, and 200 grid cells (blue curve). The combined effect of these improvements demonstrates that accuracy with hyperboloidal compactification is not restricted by the boundary treatment but by the numerical accuracy in the simulation domain.

\subsection{Three spatial dimensions}\label{sec:three}
\subsubsection{Analytical setup}\label{subsec:layer}
Consider the semilinear wave equation with a focussing power nonlinearity, $\Box_\eta u = -u^p$. We solve this equation with and without the nonlinear term by using the conformal method. The semilinear conformal wave equation reads
\be\label{eq:num} \Box_g v = \frac{1}{6} R[g] \,v -\Omega^{p-3}\, v^p,\ee
where $g$ is the conformal metric, $R[g]$ is the Ricci scalar to the metric $g$, and $v=\Omega^{-1} u$. The conformal metric $g$ is obtained from the Minkowski metric $\eta$ given in \eqref{eq:mink} by the time transformation $\tau=t-h(r)$, the spatial compactification $r=\rho/\Omega$, and the conformal rescaling $g=\Omega^2\eta$. We set the conformal factor $\Omega$ as in section \ref{sec:mult_lay}
\[ \Omega = 1- \left(\frac{\rho-R}{S-R}\right)^4 \Theta(\rho-R), \qquad  L = 1 + \frac{(\rho-R)^3 (3\rho+R)}{(S-R)^4} \Theta(\rho-R)\,.\]
This choice ensures that the coordinates $r$ and $\rho$ coincide and the Ricci scalar is continuous along the interface. The boost function is such that the outgoing characteristic speed through the layer is unity \eqref{eq:hl1om}. The conformal metric becomes 
\[ g = -\Omega^2 d\tau^2 - 2(L-\Omega^2)\,d\tau d\rho +(2L-\Omega^2)\,d\rho^2+\rho^2 d\sigma^2.\]
For $\rho<R$ we recover the Minkowski metric \eqref{eq:mink}. The characteristic speeds of spherical wave fronts are
\[ c_\pm = \frac{L\pm L - \Omega^2}{\quad 2\ L-\Omega^2}. \]
We observe that the outgoing characteristic speed $c_+$ is unity everywhere; the incoming characteristic speed $c_-$ is unity inside the domain $\rho<R$, and zero at infinity. The characteristic structure is similar to the right half of figure \ref{fig:layer_chars}. The Ricci scalar reads
\[ R[g] = \frac{6\Omega(\Omega L' - 2L\Omega')}{\rho^2 L^3}.\]
The apostrophe denotes derivative by the argument.
\subsubsection{Numerical setup}
We apply similar numerical techniques as those that have been used to test constant mean curvature foliations \cite{Zenginoglu:2010zm}. The application of hyperboloidal compactification is complicated by the requirement of a spherical grid outer boundary because spherical coordinates are ill-defined at the origin. To deal with this problem, we use the Spectral Einstein Code (SpEC) \cite{SpECWebsite}. 

\begin{figure}[ht]
\centering 
\includegraphics[width=0.6\textwidth]{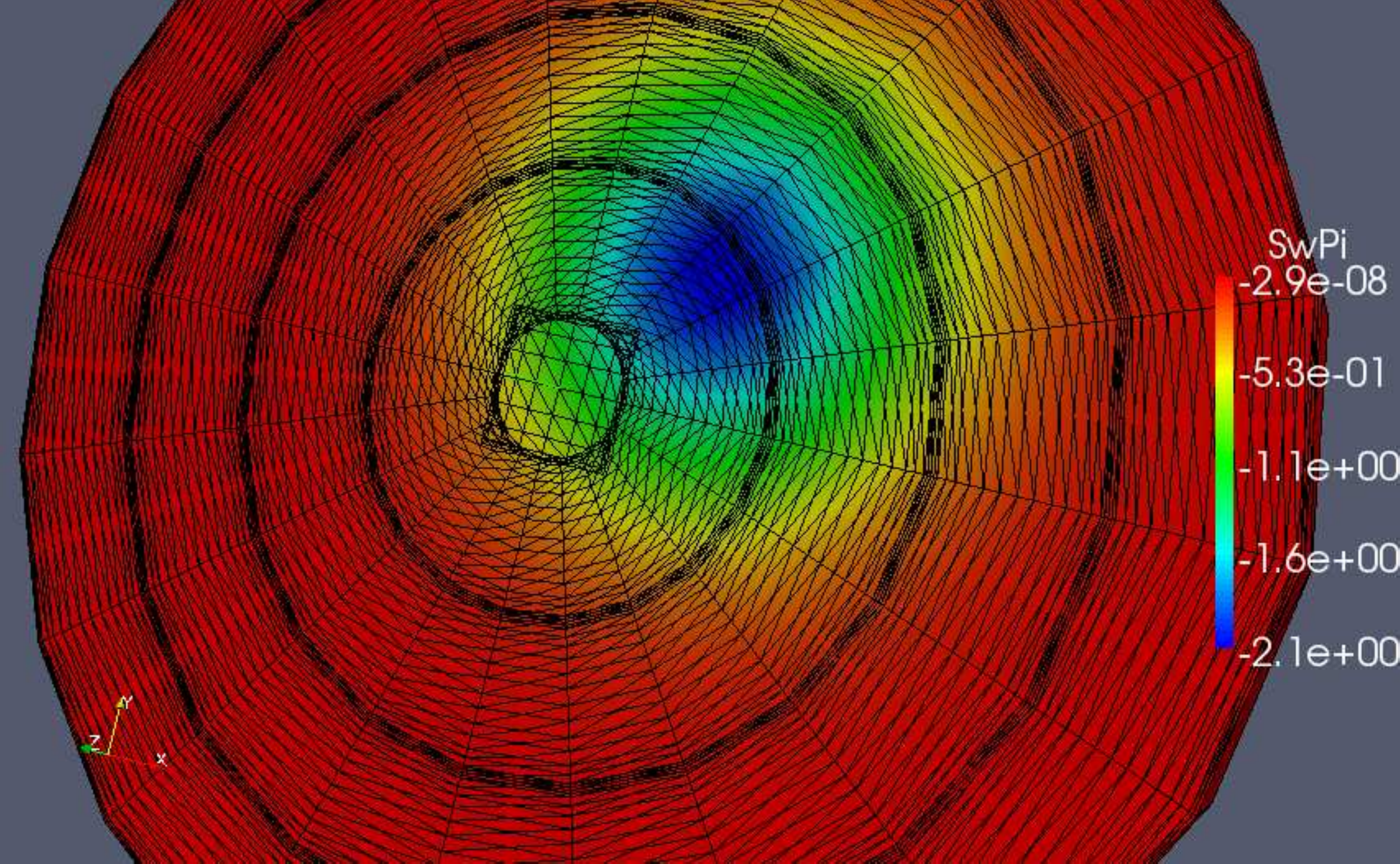}
\caption{The numerical grid for calculations in Minkowski spacetime. We have a cube around the origin on the domain $[-2,2]$ in each Cartesian direction and four spherical shells extending from $\rho=2$ to future null infinity at $\rho=20$. The colors depict off-centered Gaussian initial data for the time derivative of the scalar field. \label{fig:offc}}
\end{figure}

We cover the origin with an inner cube using Cartesian coordinates, thereby avoiding the coordinate singularity of spherical coordinates at the origin. Spherical shells starting within the cube extend to a spherical outer boundary which corresponds to null infinity (figure \ref{fig:offc}). The domain in the radial direction is $\rho\in[0,20]$. The cube around the origin has $x_i\in[-2,2]$. Four spherical shells extend from $\rho=2$ to future null infinity. 

Spatial derivatives are discretized by a pseudospectral method. For the Cartesian grid we use Chebyshev polynomials in each direction. For the shells we use Chebyshev polynomials in the radial direction and spherical harmonics in the angular directions. The  wave equation is written in first order symmetric hyperbolic form and time integration is performed with a Runge--Kutta algorithm (see \cite{Zenginoglu:2010zm} for details). The colors in figure \ref{fig:offc} depict the off-centered data prescribed for the time derivative of the scalar field. The data for the scalar field vanishes. 

The interface to the hyperboloidal layer is at $R=10$. Experiments confirm that the results reported below are not sensitive to the location of the interface if the layer is sufficiently thick to resolve an outgoing wavelength. The optimal thickness of the layer in applications will depend on certain factors such as the wavelength of the outgoing radiation or, probably, the shape of the scatterer. Future research should determine guidelines for the thickness of the layer.

\subsubsection{Results}
\begin{figure}[ht]
\center
\includegraphics[width=0.45\textwidth]{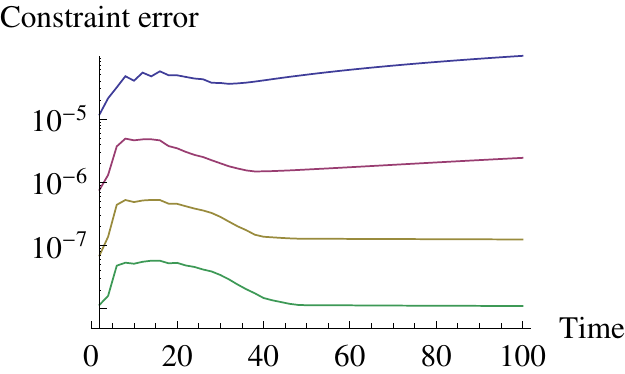}
\caption{Spectral convergence for the homogeneous wave equation with off-centered initial data. We show the $L_2$ norm of the constraint, $v_i - \partial_i v$, using $N=5,7,9,11$ collocation points (from top to bottom) in each subdomain direction.
\label{fig:conv3D}}
\end{figure}

Spatial truncation errors converge exponentially in a pseudospectral code. We show such spectral convergence in figure \ref{fig:conv3D} for a solution of the homogeneous linear wave equation with off-centered initial data. We calculate convergence in the auxiliary variable of the first order reduction $v_i:=\partial_i v$. Figure \ref{fig:conv3D} shows that the numerical difference between $v_i$ and $\partial_i v$ converges in the $L_2$ norm geometrically with $N$ collocation points in each subdomain direction. The constraint error grows slowly in time but the evolution is stable. Convergence of the constraint error is not disturbed by the presence of the boundary located at $\rho=20$. 

A stringent test for boundary methods is the inclusion of nonlinear terms in the equations. The boundary treatment for semilinear wave equations is a difficult problem, especially in three spatial dimensions \cite{Szeftel}. Therefore, we present also results for the wave equation with a cubic source term (with $p=3$ in \eqref{eq:num}). The backscatter due to the self interaction of the field violates Huygens' principle which makes the solution more interesting. Transparency boundary conditions must not eliminate all reflections from the outer boundary but only spurious ones. This makes the treatment of the artificial boundary difficult.

\begin{figure}[ht]
 \centering 
 \includegraphics[width=0.45\textwidth]{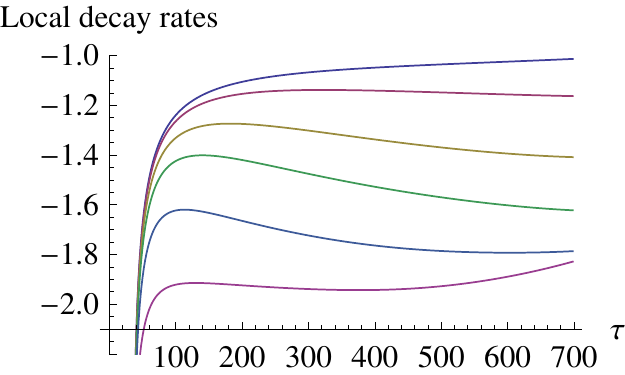}
 \caption{Local decay rates for the cubic wave equation measured by observers located from top to bottom at $\rho=\{20,19.97,19.9,19.8,19.4,17.86\}$, or equivalently at $r=\{\infty,1970,500,226,87,29\}$. The bottom curve that corresponds to the fastest decay bends up after $\tau=500$ due to accumulation of numerical errors.
\label{fig:cubic}}
\end{figure}

The backscatter by the nonlinear term leads to a late time polynomial decay of the solution: as $t^{-2}$ at a finite distance, as $t^{-1}$ at infinity. It is difficult to obtain the correct decay rates with an artificial outer boundary. The decay rate at infinity is not even accessible with standard methods. Figure \ref{fig:cubic} shows that the rates are calculated accurately with the hyperboloidal method. The local decay rate is defined as $d\ln |v(\rho, \tau)|/d\ln\tau$. The invariance of the time direction under hyperboloidal transformations implies that the rates calculated with the hyperboloidal method are equivalent to the rates calculated with the untransformed equations.

Any numerical method fails after a certain time due to accumulation of numerical errors. In figure \ref{fig:cubic} the bottom curve that corresponds to the fastest decay bends up after about $\tau=500$. This behavior is expected, and is due to accumulated numerical errors. We can delay its appearance by increasing the numerical resolution of the simulation. The long time accuracy of the solution is not limited by the boundary treatment, but solely by the accuracy of the interior calculation.
\section{Discussion}\label{sec:disc}
Hyperboloidal compactification provides a clean solution to the artificial outer boundary problem. It allows us to solve suitable hyperbolic equations accurately on unbounded domains. The underlying idea originates in studies of asymptotic structure of spacetimes in general relativity \cite{Penrose63, Friedrich83}. It is remarkable that such a practical method is available at the interface between differential geometry, conformal structure, partial differential equations, and numerical analysis. 

With the hyperboloidal layer method, in which compactification is restricted to a layer, we can apply arbitrary coordinates in an interior domain, and we gain direct quantitative access to the asymptotic solution. The numerical boundary of the simulation domain corresponds to infinity, therefore no outer boundary conditions are needed. The hyperboloidal time transformation leads to a non-vanishing coordinate speed for outgoing characteristics up to and including infinity in compactifying coordinates.

Some further advantages of the hyperboloidal method are listed as follows. Efficiency of numerical calculations improves by suitable choices of free parameters that modify the wavelength of outgoing radiation. It is straightforward to apply the method to various covariant hyperbolic equations, such as wave, Maxwell, or Yang--Mills equations. It is not necessary to calculate boundary data that typically depend on the particular equations. Including source terms and nonlinearities does not modify the boundary treatment if certain asymptotic fall-off conditions are satisfied. No separate error controlling at the outer boundary is needed because there are no errors beyond the numerical discretization. There is no overhead in software implementation of boundary routines. The method is largely independent of numerical schemes due to its geometric origin.

An important current limitation of the method is the requirement of a spherical or spheroidal grid near the outer boundary. It is not clear whether hyperboloidal compactification works with global Cartesian or cylindrical grid structures. Another limitation is source terms that lead to a singular compactification (Klein--Gordon equation, low power nonlinearities). The applicability of the method in even spatial dimensions, starting with two dimensional problems, should be studied.

The decision to apply hyperboloidal compactification can be made only after detailed comparative studies which go beyond the scope of this paper. A comparison of different boundary treatments (absorbing boundary conditions, perfectly matched layers, and hyperboloidal layers) within the same numerical setting on a simple example, such as the linear, homogeneous wave equation, would be useful. The method should be applied to Maxwell equations and its efficiency should be studied when scatterers with irregular shapes are present in the computational domain. Generalizations of the method, such as adapting the coordinate representation of the hyperboloidal layer to the shape of scatterers, might lead to very efficient computations.

The relevance of hyperboloidal compactification depends on the problem at hand. The method in its current form does not apply to non-covariant problems, however, it may be extended to deal with heterogeneous media, anisotropic elastic waves, optical waveguides, Euler equations, Navier--Stokes equations, among others.

The idea of boosting the time direction may find applications beyond compactification. It seems that transformations of the characteristic cone for hyperbolic equations are largely unexplored outside numerical relativity. The presented solution of the outer boundary problem may be just one application of such transformations. Others might be useful, for example in numerical studies of high frequency waves.

\section*{Acknowledgments} I thank Daniel Appel\"o, Eliane B\'ecache, and Fr\'ed\'eric Nataf for discussions, Larry Kidder for his help with SpEC, and  Piotr Bizo\'n, Philippe LeFloch, and Eitan Tadmor for support. This research was supported by the National Science Foundation (NSF) grant 07-07949 in Maryland, by the Marie Curie Transfer of Knowledge contract MTKD-CT-2006-042360 in Krak\'ow, by the Agence Nationale de la Recherche (ANR) grant 06-2-134423 entitled "Mathematical Methods in General Relativity" in Paris, by the NSF grant PHY-0601459 and by a Sherman Fairchild Foundation grant to Caltech in Pasadena.

\end{document}